\documentclass[11pt]{article}
\usepackage{graphicx}
\usepackage{amsmath,amsfonts,amssymb}
\usepackage{epstopdf}
\usepackage{url}
\usepackage{amsmath}
\usepackage{graphicx}
\usepackage{epsfig}
\usepackage{color}
\def\Limsup{\mathop{{\rm Lim}\,{\rm sup}}}

\def\tto{\;{\lower 1pt \hbox{$\rightarrow$}}\kern -10pt
\hbox{\raise 2pt \hbox{$\rightarrow$}}\;}
\def\Hat{\widehat}

\def\Bar{\overline}
\def\ra{\rangle}
\def\la{\langle}
\def\ve{\varepsilon}
\def\B{\Bbb B}
\def\h{\hfill\Box}
\def\R{\Bbb R}
\def\N{\Bbb N}
\def\ox{\bar{x}}
\def\ow{\bar{w}}

\def\h{\hfill\square}
\def\dn{\downarrow}
\def\O{\Omega}
\def\ph{\varphi}

\def\oR{\Bar{\R}}

\newcounter{lk}

\begin{document}

\begin{center}
\vspace*{0.3in} \textbf{GENERALIZED DIFFERENTIATION AND CHARACTERIZATIONS FOR DIFFERENTIABILITY OF INFIMAL CONVOLUTIONS}\\[2ex]
NGUYEN MAU NAM\footnote{Fariborz Maseeh Department of
Mathematics and Statistics, Portland State University, PO Box 751, Portland, OR 97207, United States (mau.nam.nguyen@pdx.edu). The research of Nguyen Mau Nam was
partially supported by the NSF under grant \#1411817 and the Simons Foundation under grant \#208785.}and  DANG VAN CUONG\footnote{Department of Mathematics, Faculty of Natural Sciences, Duy Tan University, K7/25 Quang Trung, Da Nang, Viet Nam (dvcuong@duytan.edu.vn).}
\end{center}

{\small \textbf{Abstract.} This paper is devoted to the study of generalized differentiation properties  of the infimal convolution. This class of functions covers a large spectrum of nonsmooth functions well known in the literature. The subdifferential formulas obtained unify several known results and allow us to characterize the differentiability of the infimal convolution which plays an important role in variational analysis and optimization. }

\medskip
\vspace*{0,05in} \noindent {\bf Key words.} generalized differentiation, distance function, minimal time function, infimal convolution.
\\[2ex]
{\bf Mathematical Subject Classification (2000):} 49J52, 49J53, 90C31
{\small \newtheorem{Theorem}{Theorem}[section]
\newtheorem{Proposition}[Theorem]{Proposition}
\newtheorem{Remark}[Theorem]{Remark} \newtheorem{Lemma}[Theorem]{Lemma}
\newtheorem{Corollary}[Theorem]{Corollary}
\newtheorem{Definition}[Theorem]{Definition}
\newtheorem{Example}[Theorem]{Example}
\renewcommand{\theequation}{\thesection.\arabic{equation}} }

\section{Introduction}
Throughout this paper we consider a real Banach space $X$ with a given norm $\|\cdot\|$. The dual space of $X$ is denoted by $X^*$ and the paring of an element $x^*\in X^*$ and $x\in X$ is denoted by $\la x^*,x\ra$, i.e., $\la x^*, x\ra:=x^*(x)$. The closed ball centered at $\ox$ with radius $r>0$ is denoted by $\B(\ox; r)$ and the closed unit ball of $X$ is denoted by $\B$. Given a real-valued function $\ph: X\to [0, \infty)$ and an extended-real-valued function $f: X\to \oR:=(-\infty, \infty]$ with $\mbox{\rm dom}\,f:=\{x\in X\; |\; f(x)<\infty\}\neq\emptyset$, consider the {\it infimal convolution} of $f$ and $\ph$ defined by
\begin{equation}\label{infcon}
(f\oplus \ph)(x):=\inf\{f(w)+\ph(w-x)\; |\; w\in X\}.
\end{equation}
For simplicity, we also assume that $(f\oplus \ph)(x)>-\infty$ for all $x\in X$. These are our standing assumptions throughout the paper. Under the standing assumptions, the infimal convolution (\ref{infcon}) is a real-valued function which forms an important class of nonsmooth functions containing many well-known functions in the literature. Let us emphasize its importance by some examples below.

Given a positive constant $\alpha$, consider the function $\ph(x):=\alpha\|x\|^2$. Then we obtain the \emph{quadratic infimal convolution}
\begin{equation}\label{quincon}
f_\alpha(x):=\inf\{f(w)+\alpha\|w-x\|^2\; |\; w\in X\}.
\end{equation}
The quadratic infimal convolution plays a crucial role in optimization from both theoretical and numerical aspects. It is often used to approximate a nonsmooth function by a smooth one that is convenient for applying smooth optimization schemes; see, e.g., \cite{CL1,JTZ,Nesterov} and the references therein.

The class of infimal convolutions also includes another class of functions called the \emph{minimal time function}. Let $F$ be a nonempty closed convex set that contains the origin as an interior point and let $\O$ be a nonempty subset of $X$. The minimal time function to the target set $\O$ with the dynamics $F$ is given by
\begin{equation}\label{mt}
T_F(x;\O):=\inf\{t\geq 0\; |\; (x+tF)\cap \O\neq \emptyset\}.
\end{equation}
The minimal time function (\ref{mt}) can be represented as
\begin{equation*}
T_F(x;\O)=\inf\{\rho_F(w-x)\; |\; w\in \O\}
\end{equation*}
in terms of the \emph{Minkowski function} given by $\rho_F(x):=\inf\{t\geq 0\; |\; x\in tF\}$. From this formulation we see that $T_F(x;\O)=(\delta_\O\oplus\rho_F)(x)$, where $\delta(\cdot; \O)$ is the \emph{indicator function} associated with $\O$ given by $\delta(x; \O)=0$ if $x\in \O$, and $\delta(x; \O)=\infty$ otherwise. Note that when $F$ is the closed unit ball of $X$, the minimal time function (\ref{mt}) becomes the \emph{distance function} to the set $\O$:
\begin{equation*}
d(x;\O):=\inf\{\|x-w\|\; |\; w\in \O\}.
\end{equation*}
The readers are referred to \cite{BT,BFQ,CW04a,CW04b,HNg,JH,Li,Meng,MN10a,NCA,nam,WZ,ZHJ} and the references therein for the study of the minimal time function as well as its specification to the case of the distance function.

In this paper we study generalized differentiation properties of the infimal convolution. These properties unify and provide new insights to several known results on the quadratic convolution, the minimal time function, and the distance function. We also provide new characterizations for strict differentiability of functions via generalized differentiation. Based on the results obtained, we are able to give a simple approach to study strict differentiability of the infimal convolution.

The paper is organized as follows. In Section 2 we provide some important notions and results of variational analysis used throughout the paper. General  properties of the infinal convolution are considered in Section 3. In Section 4 and Section 5 we examine generalized differentiation properties of the infimal convolution. The main attention is paid to two kinds of generalized differentiation concepts called the \emph{Fr\'echet subdifferential} and the \emph{litmiting/Mordukhovich subdifferential}. Section 6 is devoted to providing characterizations for strict differentiability  of functions and applying them to study strict differentiability of the infimal convolution.

\section{Preliminaries}

In this section we present basic notions and results of variational analysis in infinite dimensions used throughout the paper. The readers are referred to the books \cite{BZ05,CL,CL1,Mordukhovich_2006} for more details.

For a set-valued mapping
$F: X \rightrightarrows X^*$, the \textit{sequential Painlev\'e-Kuratowski upper limit} of $F$ as $x$ tends to $\bar x$ with respect to the norm topology of $X$ and the weak$^*$ topology of $X^*$ is defined by
$$
\begin{array}{rl}
\Limsup\limits_{x\rightarrow \bar x}F(x):=\Big\{x^*\in
X^*\; \big |\; & \exists \ \;x_k\rightarrow\bar x,\
x^*_k \xrightarrow{w^*} x^*,\\ &  x^*_k\in
F(x_k)\; \mbox{\rm for } k=1,2,\ldots\Big\}.
\end{array}
$$
Here  $x_k^* \xrightarrow{w^*}x^*$ means that the sequence $\{x_k^*\} \subset X^*$ converges weakly$^*$ to $x^* \in X^*$.

Given a subset $\Omega \subset X$, the notation $x \xrightarrow{\Omega}u$ means that $x \rightarrow u$ and $ x\in \Omega$. For any $x \in \Omega$ and $\varepsilon \geq 0$, the set of \textit{$\varepsilon$-normals} to $\Omega$ at $x$ is defined by
\begin{equation*}\label{varepsilon-normals}
\widehat N_\varepsilon(x;\Omega):=\Big\{ x^*\in X^*\; \big |\;  \limsup_{u \xrightarrow{\Omega}x} \dfrac{\langle x^*, u-x \rangle}{\|u-x\|} \leq \varepsilon \Big\}.
\end{equation*}
 In the case where $\varepsilon=0$, the set $\widehat N(x;\Omega):=\widehat N_0(x;\Omega)$ is called the \textit{Fr\'echet normal cone} to $\Omega$ at $x.$ If $x \not\in \Omega$, we put $\widehat N_\varepsilon(x;\Omega):=\emptyset$ for all $\varepsilon \geq 0.$

Given $\bar x \in \Omega$, the \textit{Mordukhovich normal cone} or the \textit{limiting  normal cone} to $\Omega$ at $\bar x$ is defined by
\begin{equation*}\label{1.2}
N(\bar x; \Omega) :=\Limsup_{x \rightarrow \bar x,\varepsilon \downarrow 0}\,\widehat N_\varepsilon(x; \Omega).
\end{equation*}
We also put $N(\bar x;\Omega)=\emptyset$ if $\bar x \not\in \Omega.$

Obviously, $\widehat N(x; \Omega) \subset N(x ; \Omega)$ for all $x \in \Omega$. If $\widehat N(\ox; \Omega) = N(\ox ; \Omega)$ for $\ox\in \Omega$, then one says that  $\Omega$ is {\it normally regular} at $\ox$. In the case where  $\Omega$ is a convex set, one has the following simple representation:
$$ \widehat N_\varepsilon(\bar x;\Omega)=\big\{ x^*\in X^*\mid \langle x^*, x- \bar x \rangle  \leq \varepsilon \|x-\bar x\|\; \mbox{\rm for all } x \in \Omega\big\}$$
for all $\varepsilon \geq 0$ and $\bar x \in \Omega.$ Moreover, both $\widehat N(\bar x; \Omega)$ and $N(\ox;\O)$ coincide with the convex cone to $\Omega$ at $\bar x$ in the sense of convex analysis, that is, \begin{align*}\label{normal_cone_convex_analysis} \widehat N(\bar x; \Omega)=N(\ox;\O)=\{x^*\in X^* \mid \langle x^*, x-\bar x \rangle \leq 0  \; \mbox{\rm for all } x \in \Omega\}.\end{align*}

Consider an extended-real-valued function $f: X\rightarrow \overline{\Bbb{R}}$. In the sequel, the notation $x \xrightarrow{f}\bar x$ means that $x \rightarrow \bar x$ and $f(x) \rightarrow f(\bar x).$ Given $\varepsilon\geq 0$, the $\varepsilon-$\emph{Fr\'echet subdifferential} of $f$ at $\ox\in \mbox{\rm dom}\,f$ is the set
\begin{equation*}\label{Frechetsubdifferential}
\widehat \partial_{\varepsilon} f(\bar{ x}):=\left\{x^* \in X^*\;\big |\; \liminf_{x\to \ox}\dfrac{f(x)-f(\ox)-\la x^*, x-\ox\ra}{\|x-\ox\|}\geq -\varepsilon\right\}.
\end{equation*}
The set $\widehat \partial_{0} f(\bar{ x})$ ($\varepsilon=0$) is called the \emph{Fr\'echet subdifferential} of $f$ at $\ox$ and is denoted simply by $\widehat \partial f(\bar{ x}).$

The \emph{limiting/Mordukhovich subdifferential} of $f$ at $\ox$ is defined by
\begin{align*}\label{Mordukhovichsubdifferential}
\partial f(\bar{ x}):=\Limsup_{x\xrightarrow{f}\ox, \varepsilon\dn 0}\Hat\partial_\varepsilon f(x).
\end{align*}

It follows from the definition that for any $\bar x \in \Omega$ we have
\begin{equation*}\label{Frechet_cone_indicator_function}
\widehat \partial \delta(\bar x; \Omega)=\widehat N(\bar x; \Omega)\; \mbox{\rm and }\partial \delta(\bar x; \Omega)= N(\bar x; \Omega).
\end{equation*}

The inclusion $\widehat\partial f(\bar x) \subset \partial f(\bar x)$ is valid for all $\bar x \in \mbox{\rm dom}\,f.$ If  $\widehat\partial f(\bar x)= \partial f(\bar x)$ for $\bar x\in {\rm dom}\,f$, one says that $f$ is {\it lower regular} at $\bar x$.  If $f$ is convex, then
\begin{equation*}\label{subdifferential_convex_case}\widehat\partial f(\bar x) =\partial f(\bar x)=\{x^* \in X^* \mid \langle x^*, x- \bar x \rangle \le f(x)-f(\bar x) \; \mbox{\rm for all }  x \in X\},\end{equation*}
i.e., the Fr\'{e}chet subdifferential and the Mordukhovich subdifferential of $f$ at $\bar x$ coincide with the subdifferential of $f$ at $\bar x$ in the sense of convex analysis. In particular, $f$ is lower regular at $\bar x$.

Recall that $f$ is \emph{Fr\'echet  strictly differentiable} at $\ox$ if there exists $v^*\in X^*$ such that
\begin{equation*}
\lim_{x,y\to \ox}\dfrac{f(x)-f(y)-\la v^*, x-y\ra}{\|x-y\|}=0.
\end{equation*}
The element $v^*$ is called the \emph{Fr\'echet strict derivative} of $f$ at $\ox$ and is denoted by $\nabla f(\ox)$. If $f$ is Fr\'echet strictly differentiable at $\ox$, then
\begin{equation*}
\partial f(\ox)=\Hat\partial f(\ox)=\{\nabla f(\ox)\}.
\end{equation*}

For an extended-real-valued function $f: X\to \oR$, we say that $f$ is \emph{Lipschitz continuous} on a set $D\subset\mbox{\rm dom}\,f$ with Lipschitz constant $\ell\geq 0$ if
$$|f(x)-f(w)|\leq \ell \|x-w\|\; \mbox{\rm for all }x, w\in D.$$
We also say that $f$ is \emph{locally Lipschitz continuous} at $\ox\in \mbox{\rm dom}\,f$ with constant $\ell\geq 0$ is there exists $\delta>0$ such that
$$|f(x)-f(w)|\leq \ell \|x-w\|\; \mbox{\rm for all }x,w\in \B(\ox; \delta).$$

Throughout the paper we also use other standard notations and results of variational analysis which can be found in \cite{BZ05,CL,Mordukhovich_2006}.

\section{General Properties of Infimal Convolutions}

In this section we study some general properties of the infimal convolution (\ref{infcon}). These properties will be used in the next sections.

Recall that a function $g: X\to (-\infty, \infty]$ is \emph{level bounded} if for every $\alpha\in \R$ the set
\begin{equation*}
\mathcal{L}_\alpha:=\{x\in X\; |\; g(x)\leq \alpha\}
\end{equation*}
is a bounded set in $X$. We also say that $g$ is \emph{weakly (sequentially) lower semicontinuous} on $X$ if for any $\ox\in X$ and for any sequence $\{x_k\}$ that converges weakly to $\ox$ one has
$$\liminf_{k\to \infty}g(x_k)\geq g(\ox).$$
If the weak convergence of $\{x_k\}$ is replaced by the strong convergence in the definition above, we say that $g$ is \emph{lower semicontinuous} on $X$.

\begin{Proposition} Let $X$ be a reflexive Banach space. If both $f$ and $\ph$ are weakly lower semicontinuous on $X$ and $f$ is level bounded, then $f\oplus \ph$ is weakly lower semicontinuous on $X$. In particular, it is lower semicontinuous on $X$.
\end{Proposition}
{\bf Proof.} Fix any $\ox\in X$ and any sequence $\{x_k\}$ that converges weakly to $\ox$. We will show that
\begin{equation*}
\liminf_{k \to \infty} (f\oplus \ph)(x_k)\geq (f\oplus\ph)(\ox).
\end{equation*}
Under the assumptions made, we can assume without loss of generality that $\gamma:=\liminf_{k \to \infty} (f\oplus \ph)(x_k)\in \R$ and the sequence $\{(f\oplus \ph)(x_k)\}$ converges to $\gamma$.
For every $k\in \N$, choose $w_k\in X$ such that
\begin{equation*}
f(w_k)+\ph(w_k-x_k)< (f\oplus \ph)(x_k)+1/k.
\end{equation*}
Since $\ph$ has nonnegative values and $f$ is level bounded, we see that $\{w_k\}$ is bounded in $X$, so it has a subsequence
(without relabeling) that converges weakly to $\ow\in X$. By the weak lower semicontinuity of $f$ and $\ph$,
\begin{equation*}
f(\ow)+\ph(\ow-\ox)\leq \liminf_{k\to \infty}[f(w_k)+\ph(w_k-x_k)]\leq \liminf_{k\to \infty}[(f\oplus \ph)(x_k)+1/k]=\gamma.
\end{equation*}
This implies $(f\oplus \ph)(\ox)\leq \gamma$, which completes the proof. $\h$

Recall that a function $g: X\to (-\infty, \infty]$ is called \emph{subadditive} if
$$g(x+y)\leq g(x)+g(y)\; \mbox{\rm for all }x, y\in X.$$
Given a nonempty set $D\subset X$, the function $g$ is called \emph{locally calm} at a point $\ox\in D\cap\mbox{\rm dom}\,g$ relative to $D$  if there exist constants $\ell\geq 0$  and  $\delta>0$ such that
$$|g(x)-g(\ox)|\leq \ell \|x-\ox\|\; \mbox{\rm for all }x\in \B(\ox; \delta)\cap D.$$
If the inequality above holds for all $x\in D$ instead of all $x\in \B(\ox; \delta)\cap D$, we say that $g$ is calm at $\ox$ relative to $D$. We say that $g$ is locally calm (or calm) at $\ox\in \mbox{\rm dom}\,g$ if it is locally calm (or calm) at $\ox$ relative to $X$.

Let $\ph : X \rightarrow (-\infty,\infty]$ be an extended-real-valued function. We say that $\ph$ is {\it coercive} with constant $m > 0$ on $X$ if
$$m\|x\|\leq \ph(x)\ \text{for all}\ x\in X.$$

\begin{Proposition}\label{pr3.2} Let $\ph$ be subadditive. Then
\begin{equation*}
(f\oplus\ph)(x)-(f\oplus\ph)(y)\leq \ph(y-x)
\end{equation*}
for all $x, y\in X$. Consequently, if $\ph$ is locally calm at $0$ with constant $\ell$ and $\ph(0)=0$, then $f\oplus\ph$ is locally Lipschitz continuous around any point $\ox\in X$ with Lipschitz constant $\ell$, i.e., there exists $\delta>0$ such that
\begin{equation*}
|(f\oplus\ph)(x)-(f\oplus\ph)(y)|\leq \ell\|x-y\|\; \mbox{\rm for all }x,y\in \B(\ox; \delta).
\end{equation*}
Moreover, if $\ph$ is calm at $0$ with constant $\ell$ and $\ph(0)=0$, then $f\oplus \ph$ is globally Lipschitz continuous on $X$ with constant $\ell$.
\end{Proposition}
{\bf Proof.} Fix any $x, y\in X$. Then
\begin{equation*}
f(w)+\ph(w-x)=f(w)+\ph(w-y+y-x)\leq f(w)+\ph(w-y)+\ph(y-x)\; \mbox{\rm for all }w\in X.
\end{equation*}
This implies
\begin{equation*}
(f\oplus\ph)(x)\leq f(w)+\ph(w-y)+\ph(y-x)\; \mbox{\rm for all }w\in X.
\end{equation*}
Taking the infimum with respect to $w$ on the right side yields
\begin{equation*}
(f\oplus\ph)(x)\leq (f\oplus\ph)(y)+\ph(y-x).
\end{equation*}
It follows that $(f\oplus\ph)(x)- (f\oplus\ph)(y)\leq \ph(y-x)$.

Now suppose that $\ph$ is locally calm at $0$ with constant $\ell$ and $\ph(0)=0$. Then there exists $\delta>0$ such that
\begin{equation*}
\ph(x)\leq \ell \|x\|\; \mbox{\rm for all }x\in \B(0; \delta).
\end{equation*}
For any $x,y\in \B(\ox; \delta/2),$ one has $y-x\in \B(0; \delta)$, and hence
\begin{equation*}
(f\oplus\ph)(x)- (f\oplus\ph)(y)\leq \ph(y-x)\leq \ell \|x-y\|.
\end{equation*}
This implies the locally Lipschitz continuity of $f\oplus \ph$ around $\ox$. The rest of the proof follows easily. $\h$

Let us now study the Lipschitz continuity of $f\oplus \ph$ without assuming the subadditivity of the function $\ph$.

\begin{Proposition}\label{LIP} Suppose that $f$ is bounded below on $X,$  and  $\ph$ is Lipschitz continuous and bounded above on every bounded subset of $X$. Then $f\oplus\ph$ is Lipschitz continuous on every bounded subset of $X$ under one of the following conditions:\\[1ex]
{\rm\bf (i)} $f$ is level bounded.\\
{\rm\bf (ii)} $\ph$ is level bounded.
\end{Proposition}
{\bf Proof.} Fix a bounded set $K$ and $x, y\in K$. Given any $\ox\in \mbox{\rm dom}\,f,$ one has
\begin{equation*}
(f\oplus \ph)(x)\leq f(\ox)+\ph(\ox-x)\leq f(\ox)+\sup\{\ph(u)\; |\; u\in \ox-K\}<\infty.
\end{equation*}
Define the set
\begin{equation*}
\O:=\{w\in X\; |\;\exists x\in K\;\mbox{\rm with } f(w)+\ph(w-x)<\sup_{x\in K}(f\oplus\ph)(x)+1\}.
\end{equation*}
It is not hard to see that $\O$ is nonempty and bounded under {\bf (i)} or {\bf (ii)}. For any $\varepsilon>0$ sufficiently small, choose $w\in X$ such that
\begin{equation*}
f(w)+\ph(w-x)<(f\oplus\ph)(x)+\varepsilon.
\end{equation*}
 Then $w\in \O$ and
\begin{align*}
(f\oplus\ph)(y)-(f\oplus\ph)(x)&\leq f(w)+\ph(w-y)-f(w)-\ph(w-x)+\varepsilon\\
&=\ph(w-y)-\ph(w-x)+\varepsilon\leq \ell \|x-y\|+\varepsilon,
\end{align*}
where $\ell$ is a Lipschitz constant of $\ph$ on the bounded set $\O-K$. Then we can see easily that
\begin{equation*}
|(f\oplus\ph)(y)-(f\oplus\ph)(x)|\leq \ell \|x-y\|\; \mbox{\rm for all }x,y\in K.
\end{equation*}
The proof is now complete. $\h$

For any $x\in X$, define the \emph{projection} at $x$ by
\begin{equation*}
\mathcal{P}_f^{\varphi}(x):=\{w\in X\;  |\;  f(w)+\ph(w-x)=(f\oplus\ph)(x)\}.
\end{equation*}
For simplicity, we write $\mathcal{P}(x)$ instead of $\mathcal{P}_f^{\varphi}(x)$ if no confusion occurs.

We say that $f\oplus \ph$ is\emph{ well-posed} at $\ox$ if $\mathcal{P}(\ox)$ is a singleton denoted by $\ow$ and for every sequence $\{w_k\}$ with
$$f(w_k)+\ph(w_k-\ox)\to (f\oplus\ph)(\ox),$$
we have that $\{w_k\}$ converges to $\ow$.

Following \cite{nam,ZHJ}, define the following set:
\begin{equation*}
S_0:=\{x\in X\; |\; (f\oplus \ph)(x)=f(x)\}.
\end{equation*}
The following proposition provides a sufficient condition ensuring the well-posedness of $f\oplus \ph.$
\begin{Proposition}\label{pro6} Let $\ox\in S_0.$ Assume that $f$ is calm at $\ox$ relative to $\mbox{\rm dom}\,f$ with constant $\ell,$ $\varphi$ is coercive with constant $m>\ell,$ and $\varphi(0)=0.$ Then $f\oplus\ph$ is well-posed at $\ox$.
\end{Proposition}
{\bf Proof.} Since $\ox\in S_0,$ one has $\ox\in \mathcal{P}(\ox).$ Let $\{w_k\}\subset X$ be a minimizing sequence of of $f\oplus\ph$ at $\ox$, i.e.,
$$\lim_{k\to\infty}[f(w_k)+\varphi(w_k-x)]=(f\oplus\varphi)(\ox)=f(\ox).$$
Thus, for each $\varepsilon>0$ there exists a positive integer $N$ such that if $k>N,$ then
$$f(w_k)+\varphi(w_k-\ox)\leq f(\ox)+\varepsilon\ \Leftrightarrow\ \varphi(w_k-\ox)\leq f(\ox)-f(w_k)+\varepsilon.$$
It follows that $w_k\in \mbox{\rm dom}\,f$ for such $k$, and hence
$$m\|w_k-\ox\|\leq \varphi(w_k-\ox)\leq f(\ox)-f(w_k)+\varepsilon\leq \ell\|w_k-\ox\|+\varepsilon,$$
which implies
$$\|w_k-\ox\|\leq\frac{\varepsilon}{m-\ell}.$$
Consequently, we arrive at
$$\lim_{k\to\infty}\|w_k-\ox\|=0.$$
This completes the proof of the proposition. $\h$

\section{Fr\'echet Subdifferentials of Infimal Convolutions}

In this section we develop Fr\'echet subdifferential formulas for infimal convolutions as a continuation of \cite{nam}.

\begin{Proposition} Suppose that $\ph(0)=0$, $\ph$ is coercive on $X$ with constant $m>0$ and $f$ is Lipschitz continuous on $D:=\mbox{\rm dom}\,f$ with  constant $\ell$ where $0\leq \ell <m$. Then
\begin{equation}\label{s}
S_0=\{x\in X\; |\; \mathcal{P}(x)=\{x\}\}.
\end{equation}
\end{Proposition}
{\bf Proof.} Suppose that $x\in S_0$. Then $(f\oplus\ph)(x)=f(x)=f(x)+\ph(x-x)$. It follows from the definition that $x\in\mathcal{P}(x)$. Now fix any $w\in \mathcal{P}(x)$. Then
\begin{equation*}
(f\oplus\ph)(x)=f(w)+\ph(w-x)=f(x),
\end{equation*}
which implies $m\|x-w\|\leq\ph(w-x)=f(x)-f(w)\leq \ell\|x-w\|$, so $(m-\ell)\|x-w\|=0$, which implies $x=w$. The converse also follows easily from the definition.$\h$

\begin{Example}{\rm Let $F$ be a closed bounded convex set that contains $0$ as an interior point and let $\O$ be a nonempty set. As mentioned in earlier, the minimal time function \eqref{mt} has the following representation:
\begin{equation*}
\mathcal{T}_F(x;\O)=\inf\{\rho_F(w-x)\; |\; w\in \O\}=(f\oplus\ph)(x),
\end{equation*}
 where $f(x)=\delta(x;\O)$ and $\ph(x)=\rho_F(x)$. Then $\rho_F(0)=0$ and $\rho_F(\cdot)$ is coercive with constant $m:=\|F\|^{-1}$, where $$\|F\|:=\sup\{\|f\|\; |\; f\in F\}.$$
Moreover, $f$ is Lipschitz continuous on $D:=\mbox{\rm dom}\,f$ with constant $\ell=0$. It is not hard to see that $S_0=\O$.}
\end{Example}

Let us present below a result on Fr\'echet-type subdifferential for the infimal convolution (\ref{infcon}) obtained in \cite{nam}.
\begin{Theorem}\label{t1} Suppose that $\ph(0)=0$ and consider the set $S_0$ given by {\rm (\ref{s})} with $\ox\in S_0$. \\[1ex]
{\rm\bf(i)} Given $\varepsilon\geq 0$, one has
\begin{equation*}
\Hat\partial_\varepsilon (f\oplus\ph)(\ox)\subset \Hat\partial_\varepsilon f(\ox)\cap \big[-\Hat\partial_\varepsilon \ph(0)\big].
\end{equation*}
{\rm\bf(ii)} Suppose that and $\ph$ is coercive on $X$ with constant $m>0$ and $f$ is calm at $\ox$ relative to $D:=\mbox{\rm dom}\,f$ with constant $\ell$ where $0\leq \ell <m$. Given $\varepsilon\geq 0$ and $x^*\in \Hat\partial_\varepsilon f(\ox)\cap \big[-\Hat\partial_\varepsilon \ph(0)\big]$, one has
\begin{equation*}
x^*\in \Hat\partial_{\alpha\varepsilon}(f\oplus\ph)(\ox), \mbox{\rm where }\alpha:=2(\|x^*\|+m)(m-\ell)^{-1}+1.
\end{equation*}
Moreover,
\begin{equation*}
\Hat\partial (f\oplus\ph)(\ox)= \Hat\partial f(\ox)\cap \big[-\Hat\partial \ph(0)\big].
\end{equation*}
\end{Theorem}

Now we consider the case where the reference point is not necessarily in the set $S_0$.
\begin{Proposition}\label{fr}  Given $\ox\in X$, suppose that $\mathcal{P}(\ox)$ is nonempty. Then
$$
\widehat{\partial_{\varepsilon}} (f\oplus \ph)(\ox )\subset \bigcap_{w\in \mathcal{P}(\ox)} \left(\widehat{\partial_{\varepsilon}} f(w)\cap [-\Hat\partial_\varepsilon \ph(w-\ox)]\right).$$
\end{Proposition}
{\bf Proof.} Fix any $x^*\in \Hat\partial_{\varepsilon}(f\oplus \ph)(\ox)$ and $w\in \mathcal{P}(\ox)$. Then for any $\eta>0$ there exists $\delta>0$ such that
\begin{align}\label{subdef}
\la x^*, x-\ox\ra &\leq (f\oplus \ph)(x)-(f\oplus \ph)(\ox)+(\ve+\eta)\|x-\ox\|\; \mbox{\rm whenever }\|x-\ox\|<\delta.
\end{align}
Fix any $z\in X$ with $\|z-w\|<\delta$. Then $\|z-w+\ox-\ox\|<\delta$, and hence we can apply \eqref{subdef} with $x$ replaced by $z-w+\ox$ to obtain
\begin{align*}
\la x^*, z-w\ra &\leq (f\oplus \ph)(z-w+\ox)-f(w)-\ph(w-\ox)+(\ve+\eta)\|z-w\|\\
&\leq f(z)+\ph(w-\ox)-f(w)-\ph(w-\ox)+(\ve+\eta)\|z-w\|\\
&=f(z)-f(w)+(\ve+\eta)\|z-w\|.
\end{align*}
It follows that $x^*\in \Hat\partial_{\varepsilon} f(w)$.

Moreover, from \eqref{subdef}, one has
\begin{align*}
\la x^*, x-\ox\ra &\leq (f\oplus \ph)(x)-(f\oplus \ph)(\ox)+(\ve+\eta)\|x-\ox\|\\
&=(f\oplus \ph)(x)-f(w)-\ph(w-\ox)+(\ve+\eta)\|x-\ox\|\\
&\leq f(w)+\ph(w-x)-f(w)-\ph(w-\ox)+(\ve+\eta)\|x-\ox\|\\
&=\ph(w-x)-\ph(w-\ox)+(\ve+\eta)\|x-\ox\|\; \mbox{\rm whenever }\|x-\ox\|<\delta.
\end{align*}
It follows that $-x^*\in \Hat\partial_\ve \ph(w-\ox)$. The proof is now complete. $\h$

Let us now consider the case where $\ph$ is subadditive and positively homogeneous.

\begin{Proposition}\label{pro2.3}
Suppose that $\varphi$ is subadditive and positively homogeneous. Let $\ox\in X$ and $\ow\in \mathcal{P}(\ox). $ Then, for each $t\in(0,1],$ we have $\ow\in \mathcal{P}(t\ow+(1-t)\ox).$  Consequently,
$$
(f\oplus\ph)(t\ow+(1-t)\ox)=(1-t)(f\oplus\ph)(\ox)+tf(\ow)\ \text{for each}\ t\in(0,1].
$$
\end{Proposition}
{\bf Proof.} Let $t\in(0,1]$ and set $x_t:=t\ow+(1-t)\ox.$ Since $\varphi$ is subadditive and positively homogeneous, for every $w\in X,$ we have
$$\begin{aligned}
f(\ow)+\varphi(\ow -x_t)&=f(\ow)+\varphi[(1-t)(\ow-\ox)]\\
&=f(\ow)+(1-t)\varphi(\ow-\ox)\\
&=f(\ow )+\varphi(\ow-\ox  )-t\varphi(\ow -\ox)\\
&=(f\oplus\ph)(\ox )-t\varphi(\ow-\ox)\\
&\leq f(w)+ \varphi(w-\ox )-t\varphi(\ow-\ox)\\
&=f(w )+\varphi(w-\ox )-\varphi(x_t-\ox )\\
&\leq f(w)+\varphi(w-x_t).
\end{aligned} $$
It means that $\ow\in \mathcal{P}(x_t).$ Consequently,
$$\begin{aligned}
(f\oplus\ph)(x_t)&=f(\ow )+\varphi(\ow-x_t )=f(\ow )+\varphi[(1-t)(\ow-\ox  )]\\& =(1-t)f(\ow )+(1-t)\varphi(\ow-\ox)+tf(\ow)\\&=(1-t)(f\oplus\ph)(\ox )+tf(\ow).
\end{aligned} $$
The proof is now complete. $\h$

\begin{Theorem}\label{thr2.4} Suppose that $\varphi$ is subadditive and positively homogeneous. Let $\ox\in X$ satisfy $\mathcal{P}(\ox)\ne \emptyset.$ Then we have
$$\widehat{\partial_{\varepsilon}} (f\oplus \ph)(\ox )\subset \bigcap_{w\in \mathcal{P}(\ox)}\bigcap_{t\in(0,1]}\left([\widehat{\partial_{\varepsilon}} (f\oplus \ph)(tw+(1-t)\ox )]\cap [-\Hat\partial_\varepsilon \ph(w-\ox)]\right).$$
\end{Theorem}
{\bf Proof.} Let $w\in \mathcal P(\ox)$ and let $t\in(0,1].$ We will show that
\begin{equation*}
\widehat{\partial_{\varepsilon}} (f\oplus \ph)(\ox )\subset \widehat{\partial_{\varepsilon}} (f\oplus \ph)(tw+(1-t)\ox ).
\end{equation*}
Fix any $x^*\in\widehat{\partial_{\varepsilon}}(f\oplus \ph)(\ox)$ and let $\eta>0.$ Then there exists $\delta>0$ such that
\begin{equation}\label{eq1}
\langle x^*,x-\ox \rangle\leq (f\oplus \ph)(x)-(f\oplus \ph)(\ox )+(\varepsilon+\eta)\|x-\ox \|\; \mbox{\rm for all } x\in\B(\ox,\delta).
\end{equation}
Let $x_t:=tw+(1-t)\ox.$ For any $u\in \B(x_t,\delta),$ we have $u-x_t=u-t(w-\ox )-\ox\in \delta\B$, and so $u-t(w-\ox )\in  \B(\ox,\delta).$ Applying (\ref{eq1}) with $x:= u-t(w-\ox )$ yields
\begin{equation*}\label{eqthr2.4.6} \langle x^*,u-x_t\rangle\leq (f\oplus \ph)(u-t(w-\ox ))-(f\oplus \ph)(\ox )+(\varepsilon+\eta)\|u-x_t\|. \end{equation*}
Since $\varphi$ is subadditive and positively homogeneous, Proposition \ref{pr3.2} implies that
\begin{equation*}
(f\oplus \ph)(u-t(w-\ox ))\leq (f\oplus \ph)(u)+t\varphi(w-\ox).
\end{equation*}
It follows that
\begin{equation}\label{eq2}
\langle x^*,u-x_t\rangle \leq (f\oplus \ph)(u)+t\varphi(w-\ox)-(f\oplus \ph)(\ox)+(\varepsilon+\eta)\|u-x_t\|.
\end{equation}
By Proposition \ref{pro2.3},
\begin{align*}
(f\oplus \ph)(x_t)&=(1-t)(f\oplus \ph)(\ox )+tf(w)\\
&=(f\oplus \ph)(\ox)-t[(f\oplus \ph)(\ox)-f(w)]\\
&=(f\oplus \ph)(\ox )-t\varphi(w-\ox ).
\end{align*}
Substituting into (\ref{eq2}) yields
$$
\langle x^*,u-x_t\rangle \leq (f\oplus \ph)(u)-(f\oplus \ph)(x_t)+(\varepsilon+\eta)\|u-x_t\|.
$$
This implies $x^*\in\widehat{\partial_{\varepsilon}} (f\oplus \ph)(x_t)$.

It follows from Proposition \ref{fr} that $-x^*\in \Hat\partial_\varepsilon\ph(w-\ox)$ and we have justified the theorem.  $\h$

\section{Limiting Subdifferentials of Infimal Convolutions}

Given $\ox\in X$ and $\eta>0$, define
\begin{equation*}
\mathcal{P}(\ox; \eta):=\{w\in X\; |\; f(w)+\ph(w-x)<(f\oplus\ph)(\ox)+\eta\}.
\end{equation*}
Note that this set is always nonempty.

\begin{Lemma}\label{lmdomfS0} Suppose that $\ph(0)=0$, $\ph$ is coercive on $X$ with constant $m>0$, and $f$ is Lipschitz continuous on $D:=\mbox{\rm dom}\,f$ with constant $\ell,$ where $0\leq \ell <m$. Then $\mbox{\rm dom}\,f\subset S_0$. In particular,
\begin{equation*}
\mathcal{P}(\ox; \eta)\subset S_0.
\end{equation*}
\end{Lemma}
{\bf Proof.} Fix any $x\in \mbox{\rm dom}\,f$.  If, by contradiction, $x\notin S_0$, then
\begin{equation*}
(f\oplus \ph)(x)<f(x).
\end{equation*}
Then there exists $w\in X$ such that $f(w)+\ph(w-x)<f(x)$, and hence $\ph(w-x)<f(x)-f(w)\leq\ell\|x-w\|$. It follows that
$$m\|w-x\|<\ell\|w-x\|,$$
So $(m-\ell)\|w-x\|<0$. This is a contradiction. $\h$\\

We recall the well-known Ekeland variational principle; see, e.g., \cite{Figueiredo}.
\begin{Proposition}[Ekeland's variational principle]
  Let $(E,d)$ be a complete metric space and let $\phi:E\to\oR$ be a proper lower semicontinuous function that is bounded below. Let $\widetilde{\eta}>0$ and $\widetilde{w}\in E$ such that
  \begin{equation}\label{eqpro4.1}\phi(\widetilde{w})\leq \inf_{w\in E}\phi(w)+\widetilde{\eta}.\end{equation}
  Then for any $\lambda>0$ there exists $\ow\in E $ satisfying
 $$ \phi(\ow )\leq\phi(\widetilde{w}),\ d(\ow,\widetilde{w})\leq \lambda$$
  and
$$\phi(\ow )\leq \phi(w)+\frac{\widetilde{\eta}}{\lambda}d(w,\ow) \text{ for all } w\in E.$$
\end{Proposition}
\begin{Lemma}\label{prop5.1} Suppose that $\ph$ is lower semicontinuous. Let  $\varepsilon>0, \eta>0,\ox\in X,$ and $x^*\in \Hat\partial_\varepsilon(f\oplus\ph)(\ox).$
Then there exist  $\widetilde w\in \mathcal P(\ox,\eta^2)$ and $\ow\in X $ such that
\begin{equation*}\|\ow-\widetilde w\|<\eta \ \text{and}\ \ x^*\in -\Hat\partial_{\varepsilon+\eta}\ph(\ow-\ox).\end{equation*}
\end{Lemma}
{\bf Proof.} Fix any $x^*\in \Hat\partial_\varepsilon(f\oplus\ph)(\ox).$
It follows from the definition of $\Hat\partial_{\varepsilon}(f\oplus\ph)(\ox) $ that  there exists $0<\delta<\frac{\eta}{2} $ such that
\begin{equation*}
\la x^*, x-\ox\ra \leq (f\oplus\ph)(x)-(f\oplus\ph)(\ox)+(\varepsilon+\frac{\eta}{2} ) \|x-\ox\|\; \mbox{\rm for all }x\in \B(\ox,\delta).
\end{equation*}
Let $0<\widetilde\eta<\frac{\delta}{2}.$ Fix $\widetilde w\in X$ such that
$$f(\widetilde w)+\ph(\widetilde w-\ox)<(f\oplus \ph)(\ox)+\widetilde\eta^2<(f\oplus \ph)(\ox)+\eta^2.$$
For any $w\in \B(\widetilde w,\delta),$ one has $\widetilde w-w+\ox\in\B(\ox,\delta).$ Therefore,
\begin{align*}
\la x^*, \widetilde w-w\ra &\leq (f\oplus\ph)(\widetilde w-w+\ox)-(f\oplus\ph)(\ox)+(\varepsilon+\frac{\eta}{2} ) \|\widetilde w-w\|\\
&\leq (f\oplus \ph)(\widetilde w-w+\ox)-f(\widetilde w)-\ph(\widetilde w-\ox)+\widetilde\eta^2+(\varepsilon+\frac{\eta}{2} ) \|\widetilde w-w\|\\
&\leq f(\widetilde w)+\ph(w-\ox)-f(\widetilde w)-\ph(\widetilde w-\ox)+\widetilde\eta^2+(\varepsilon+\frac{\eta}{2} )\|\widetilde w-w\|\\
&= \ph(w-\ox)-\ph(\widetilde w-\ox)+\widetilde\eta^2+(\varepsilon +\frac{\eta}{2} )\|\widetilde w-w\|.
\end{align*}
Define $\phi(w):=-\la x^*,\widetilde w-w\ra +\ph(w-\ox)-\ph(\widetilde w-\ox)+\widetilde\eta^2+(\varepsilon+\frac{\eta}{2})\|\widetilde w-w\|,$ where $w\in\B(\widetilde w,\delta).$

It is easy to show that $\phi$ is lower semicontinuous, $\phi(\widetilde w)=\widetilde\eta^2,$ and $\phi(w)\geq 0$ for all $w\in \B(\widetilde w,\delta)$.
By the Ekeland variational principle applied to $\phi$ on $\B(\widetilde w,\delta)$, there exists $\ow\in \B(\widetilde w, \delta)$ such that
\begin{equation}\label{eqpr5.12} \|\widetilde w-\ow\|<\widetilde\eta<\eta
\end{equation}
and
\begin{equation}\label{eqpr5.13} \phi(\ow)\leq \phi(w)+\widetilde\eta\|w-\ow\| \text{ for all } w\in\B(\widetilde w,\delta).
\end{equation}
By (\ref{eqpr5.13}), we have
\begin{equation}\label{eqpr5.14}\begin{aligned}-\la x^*,w-\ow \ra&\leq \ph(w-\ox)-\ph(\ow-\ox)+(\varepsilon +\frac{\eta}{2} )\|w-\ow\|+\widetilde\eta\|w-\ow\|\\
 &\leq \ph(w-\ox)-\ph(\ow-\ox)+(\varepsilon +\eta)\|w-\ow\| \text{ for all } w\in \B(\widetilde w,\delta).\end{aligned} \end{equation}
 Since  $0<\widetilde\eta<\frac{\delta}{2}$ and using (\ref{eqpr5.12}),  for any $w\in\B(\ow,\widetilde\eta)$ one has
 $$\|w-\widetilde w\|\leq\|w-\ow\|+\|\widetilde w-\ow\|\leq \widetilde\eta+\widetilde\eta<\delta.$$
 It follows that $\B(\ow,\widetilde\eta)\subset\B(\widetilde w,\delta).$ Thus, (\ref{eqpr5.14}) holds for all $w\in\B(\ow,\widetilde\eta)$ and so $x^*\in-\Hat\partial_{\varepsilon+\eta}\ph(\ow-\ox).$  $\h$

\begin{Lemma}\label{lm5} Suppose that $f$ is lower semicontinuous.
Let $\varepsilon>0, \eta>0,$ $\ox\in X,$ and $x^*\in\widehat{\partial_{\varepsilon}}(f\oplus\varphi)(\ox). $
Then there exist $\widetilde{w}, \ow\in \mbox{\rm dom}\,f$  such that
$$\|\ow-\widetilde{w} \|\leq\eta,\ x^*\in\widehat{\partial}_{\varepsilon+\eta}f(\ow).$$
 If we assume   further that $\varphi$ is subadditive, then
\begin{equation}\label{eqlm5.2}f(\widetilde w)+\varphi(\ow-\ox )\leq (f\oplus\varphi)(\ox )+\ph(\ow-\widetilde w)+\eta.\end{equation}
\end{Lemma}
{\bf Proof.} Since $x^*\in\widehat{\partial_{\varepsilon}}(f\oplus\varphi)(\ox)$, given any $\eta>0$, there exists $\delta>0$ such that
\begin{equation}\label{eqlm5.3}
\langle x^*,x-\ox\rangle\leq (f\oplus\varphi)(x)-(f\oplus\varphi)(\ox )+\left(\varepsilon+\frac{\eta}{2}\right)\|x-\ox\|\ \text{for all}\ x\in\mathbb B(\ox,\delta).
\end{equation}
Set $\widetilde{\eta}:=\min\{\frac{\eta}{2},\frac{\delta}{2},1\}$ and choose $\widetilde{w}\in X$ such that
\begin{equation}\label{eqlm5.4}f(\widetilde{w})+\varphi(\widetilde{w}-\ox)\leq (f\oplus\varphi)(\ox)+\widetilde{\eta}^2. \end{equation}
This implies $\widetilde w\in\mathcal P(\ox,\eta)\subset\mbox{\rm dom}\,f.$ Now we consider the metric space $\mathbb B(\widetilde{w},\delta)$ and the function $\phi:\mathbb B(\widetilde{w},\delta)\to \overline{\mathbb R} $ defined by
$$
\phi(w):=-\langle x^*,w-\widetilde{w}\rangle +f(w)-f(\widetilde{w})+\widetilde{\eta}^2+\left(\varepsilon+\frac{\eta}{2}\right)\|w-\widetilde{w}\|.
$$
Obviously, $\mathbb B(\widetilde{w},\delta)$ is a complete metric space and $\phi$ is a lower semicontinuous function. Observe that $\phi(\widetilde w)=\widetilde{\eta}^2.$ Fix any $w\in \mathbb B(\widetilde{w},\delta).$ Then $w-\widetilde w+\ox\in\mathbb B(\ox,\delta).$ It follows from (\ref{eqlm5.3}) and (\ref{eqlm5.4}) that
$$
\begin{aligned}
\langle x^*,w-\widetilde w\rangle &\leq (f\oplus\varphi)(w-\widetilde w+\ox)-(f\oplus\varphi)(\ox )+\left(\varepsilon +\frac{\eta}{2}\right)\|w-\widetilde w\|\\
&\leq (f\oplus\varphi)(w-\widetilde w+\ox)-f(\widetilde w)-\varphi(\widetilde w-\ox)+\widetilde{\eta}^2+
\left(\varepsilon +\frac{\eta}{2}\right)\|w-\widetilde w\|\\
&\leq f(w)+\varphi(\widetilde w-\ox)-f(\widetilde w)-\varphi(\widetilde w-\ox)+\widetilde{\eta}^2+
\left(\varepsilon +\frac{\eta}{2}\right)\|w-\widetilde w\|\\
&=f(w)-f(\widetilde w)+\widetilde{\eta}^2+ \left(\varepsilon +\frac{\eta}{2}\right)\|w-\widetilde w\|.
\end{aligned}
$$
Hence $\phi(w)\geq 0$ on $\mathbb B(\widetilde{w},\delta)$. Similar to the proof of Lemma \ref{prop5.1}, we can apply the Ekeland variational principle and find $\ow \in \mathbb B(\widetilde{w},\delta)$ such that
$$ \|\widetilde w-\ow \|\leq \widetilde{\eta}\leq \eta$$
 and
 \begin{equation}\label{eqlm5.7}\phi(\ow )\leq \phi(w)+\widetilde{\eta}\|w-\ow\| \text{ for all } w\in \mathbb B(\widetilde w,\delta).\end{equation}
 By the construction of $\phi(w),$ (\ref{eqlm5.7}) implies that  $\ow\in{\rm dom}f$  and
 \begin{equation}\label{eqlm5.8}
 \langle x^*,w-\ow \rangle \leq f(w)-f(\ow )+(\varepsilon+\eta)\|w-\ow \| \text{ for all }
 w\in \mathbb B(\widetilde w,\delta).
 \end{equation}
 Since
 $$\|w-\widetilde w\|\leq \|w-\ow \|+\|\ow-\widetilde w \|\leq 2\widetilde{\eta}\leq \delta \text{ for all } w\in \mathbb B(\ow ,\widetilde{\eta}),$$
 one has $\mathbb B(\ow,\widetilde{\eta})\subset \mathbb B(\widetilde w,\delta).$ This, together with (\ref{eqlm5.7}), implies that (\ref{eqlm5.8}) holds for all $w\in \mathbb B(\ow  ,\widetilde{\eta})$ and so $x^*\in \widehat{\partial}_{\varepsilon+\eta}f(\ow ).$

 If $\varphi$ is subadditive, it follows from (\ref{eqlm5.4}) that
 $$\begin{aligned}
 f(\widetilde w)+\varphi(\ow-\ox )&\leq f(\widetilde w)+\varphi(\ow-\widetilde w )+\varphi(\widetilde w-\ox)\\
 &\leq (f\oplus\varphi)(\ox)+\widetilde{\eta}^2+\varphi(\ow-\widetilde w)\\
  &\leq (f\oplus\varphi)(\ox)+\ph(\ow-\widetilde w)+\eta.
 \end{aligned} $$
 Hence (\ref{eqlm5.2}) holds and the proof is complete.$\h$

 \begin{Theorem}\label{thr8} Let $\ox\in S_0.$ Suppose that $\ph$ is coercive on $X$ with constant $m>0$ and $f$ is a lower  semicontinuous function on $X$
 which is Lipschitz continuous on $D:=\mbox{\rm dom}\,f$ with constant $\ell$ where $0\leq \ell <m$. Suppose further that $\ph$ is subadditive and continuous at $0$ with $\ph(0)=0$. Then we have
\begin{equation}\label{eqthr8.1}\partial (f\oplus\varphi)(\ox)\subset \partial f(\ox)\bigcap [-\partial \varphi(0)]. \end{equation}
Moreover,
\begin{equation}\label{eqthr8.2}\partial (f\oplus\varphi)(\ox)= \partial f(\ox)\bigcap [-\partial \varphi(0)] \end{equation}
if we assume additionally that $\ph$ positively homogeneous and  one of the following conditions holds:
\begin{enumerate}
\item[\rm\bf (i)] $X$ is  finite dimensional.
\item[\rm\bf (ii)] $f$ is lower regular at $\ox$.
\end{enumerate}
\end{Theorem}
{\bf Proof.} Let $x^*\in\partial (f\oplus\varphi)(\ox).$ Then there exist sequences $\varepsilon_k\dn 0$, $\{x_k\}\subset X,\{x_k^*\}\subset X^*$ such that $x_k\xrightarrow{f\oplus\ph}\ox$, $x^*_k\xrightarrow{w^*}x^*$ and $x^*_k\in \Hat\partial_{\varepsilon_k} (f\oplus\ph)(x_k)$.
We will first show that $x^*\in \partial f(\ox).$ By Lemma \ref{lm5}, there exist $\ow_k, \widetilde w_k\in \mbox{\rm dom}\,f$ such that
$$
\|\ow_k-\widetilde w_k\|\leq \frac{1}{k},\ x_k^*\in \widehat{\partial}_{\varepsilon_k+\frac{1}{k}}f(\ow_k),
$$
and
$$f(\widetilde w_k)+\varphi(\ow_k-x_k )\leq (f\oplus\varphi)(x_k )+\ph(\ow_k-\widetilde w_k)+1/k.$$
This implies
\begin{equation*}
m\|x_k-\ow_k\|\leq (f\oplus\varphi)(x_k )-f(\widetilde w_k)+\ph(\ow_k-\widetilde w_k)+1/k.
\end{equation*}
It follows that
\begin{align*}
m\limsup \|x_k-\ow_k\|&\leq \limsup [(f\oplus\varphi)(x_k )-f(\widetilde w_k)+\ph(\ow_k-\widetilde w_k)+1/k]\\
&\leq \limsup [(f\oplus\varphi)(x_k)-f(\widetilde w_k)]\\
&\leq \limsup [f(\ox)-f(\widetilde w_k)]\\
&\leq \ell \limsup \|\ox -\widetilde w_k\|\leq \ell \limsup (\|\ox-x_k\|+\|x_k -\ow_k\|+\|\ow_k-\widetilde w_k\|)\\
&\leq \ell \limsup \|x_k-\ow_k\|.
\end{align*}
Thus, $\limsup \|x_k-\ow_k\|=0$, and hence $\ow_k\to \ox$ as $k\to \infty$.
 Since both $\ow_k$ and $\ox$ are in $\mbox{\rm dom}\,f$,
$$|f(\ow_k)-f(\ox)\|\leq \ell \|\ow_k-\ox\|\to 0.$$
Therefore, $x^*\in \partial f(\ox)$.

Let us now show that $x^*\in -\partial \ph(0)$. By Lemma \ref{prop5.1}, there exist $\widetilde{w}_k\in X,$ $\ow_k\in X $ such that
$$f(\widetilde{w}_k)+\ph(\widetilde{w}_k-x_k)<(f\oplus\ph)(x_k)+1/k^2,\ \|\widetilde{w}_k-\-\ow_k \|<\frac{1}{k}, x_k^*\in\Hat\partial_{\varepsilon+\frac{1}{k} } \ph(\ow_k-x_k).$$

Similar to the proof above, we can show that $\widetilde w_k\to \ox$, and hence $\ow_k\to \ox$.
Then $\ph(\ow_k-x_k)\to \ph(0)$ by the continuity of $\ph$ at $0$, and hence $x^*\in -\partial\ph(0)$. T
herefore, $x^*\in\partial f(\ox)\cap[-\partial \ph(0)]$ and (\ref{eqthr8.1}) has been proved.

To prove (\ref{eqthr8.2}), it suffices to show that
$$ \partial f(\ox)\bigcap [-\partial \varphi(0)]\subset\partial (f\oplus\varphi)(\ox).$$
Let $x^*\in\partial f(\ox)\bigcap[-\partial\ph(0)].$ Then there exist $\varepsilon_k\in[0,1],$ $x_k\in X,$ and $x_k^*\in X^*$
such that
$$\varepsilon_k\downarrow0,\ x_k\xrightarrow{f}\ox,\ x_k^*\xrightarrow{w^*} x^*\ \text{and}\ x_k^*\in\Hat\partial_{\varepsilon_k}f(x_k).$$
Since $x_k\xrightarrow{f}\ox,$ for any $\varepsilon>0$ there exists $k_1>0$ such that $|f(x_k)-f(\ox)|\leq \varepsilon$ for all $k>k_1.$
By Lemma \ref{lmdomfS0}, $x_k\in \mbox{\rm dom}\,f\subset S_0$  for such $k$ and hence  $x_k\xrightarrow{f\oplus\ph}\ox.$

Using property {\rm\bf(i)}, set $\sigma_k:=\|x_k^*-x^*\|.$ Since $X$ is  finite dimensional  and $\ph$ is convex,
$$\begin{aligned}\langle -x_k^*,x\rangle&= \langle -x^*,x\rangle+ \langle -x_k^*+x^*,x\rangle\\
&\leq \ph(x)+\langle -x_k^*+x^*,x\rangle\leq \ph(x)+\sigma_k\|x\| \text{ for all } x\in X. \end{aligned}$$
This implies $-x_k^*\in \Hat\partial_{\sigma_k}\ph(0).$ Set $\delta_k:=\max\{\varepsilon_k,\sigma_k\}.$
Then $x_k^*\in \Hat\partial_{\delta_k}f(x_k)\bigcap[-\Hat\partial_{\delta_k}\ph(0)],$ and $\delta_k\downarrow0.$
Since  $x_k\in S_0$, it follows from Theorem \ref{t1} that $x_k^*\in \Hat\partial_{\alpha_k\delta_k}(f\oplus\ph)(x_k),$ where
\begin{equation}\label{eqthr8.6}\alpha_k:=2(\|x_k^*\|+m)(m-\ell)^{-1}+1.\end{equation}
Taking into account that $\{x_k^*\}$ is bounded,  (\ref{eqthr8.6}) shows that $\eta_k:=\alpha_k\delta_k\downarrow0.$   So $$\eta_k\downarrow0,\ x_k\xrightarrow{f\oplus\ph}\ox,\ x_k^*\xrightarrow{w^*}x^*\ \text{with}\ x_k^*\in\Hat\partial_{\eta_k}(f\oplus \ph)(x_k).$$
It follows that $x^*\in\partial(f\oplus\ph)(\ox),$ and (\ref{eqthr8.2}) has been proved.

Now we assume that  {\rm\bf(ii)} holds. Since $f$ is lower regular and $\ph$ is convex,
$$\partial f(\ox)\bigcap[-\partial \ph(0)]=\Hat\partial f(\ox)\bigcap [-\Hat\partial\ph(0)].$$
It is followed from Theorem \ref{t1} that
$$x^*\in \partial f(\ox)\bigcap[-\partial \ph(0)]=\Hat\partial f(\ox)\bigcap [-\Hat\partial\ph(0)]=\Hat\partial (f\oplus\ph)(\ox)\subset \partial(f\oplus\ph)(\ox).$$
The proof is complete. $\h$

Let us now focus on the case where the reference point is not necessarily in the set $S_0$.
\begin{Definition} {\rm The mapping $\mathcal{P}$ is said to be \emph{inner semicompact} at $\ox$ if $\mathcal{P}(\ox)\neq\emptyset$ and for every sequence $\{x_k\}\subset X$ converging to $\ox$, there is a sequence $\{w_k\}$ with each $w_k\in \mathcal{P}(x_k)$ that contains a subsequence converging to $\ow\in \mathcal P(\ox)$.}
\end{Definition}

\begin{Proposition}\label{lmt3} Suppose that $\ph$ is continuous at $w-\ox$ for every $w\in \mathcal{P}(\ox)$    and $\mathcal{P}$ is inner semicompact at $\ox$. Then
\begin{equation*}
\partial(f\oplus \ph)(\ox)\subset \bigcup_{\ow\in \mathcal{P}(\ox)}\left(\partial f(\ow)\cap [-\partial \ph(\ow-\ox)]\right).
\end{equation*}
\end{Proposition}
{\bf Proof.} Fix any $x^*\in \partial(f\oplus \ph)(\ox)$. Then there exist sequences $x_k\xrightarrow{f\oplus\ph}\ox$, $\varepsilon_k\dn 0$, $x^*_k\xrightarrow{w^*}x^*$ with $x^*_k\in \Hat\partial_{\varepsilon_k}(f\oplus\ph)(x_k)$. Then there exists a sequence $\{w_k\}$ with $w_k\in \mathcal{P}(x_k)$ that contains a subsequence (without relabeling) converging to $\ow\in \mathcal{P}(\ox)$. By Proposition \ref{fr},
\begin{equation*}
x^*_k\in \Hat\partial_{\varepsilon_k} f(w_k)\cap [-\Hat\partial_{\varepsilon_k}\ph(w_k-x_k)].
\end{equation*}
Since $$f(w_k)+\ph(w_k-x_k)=(f\oplus\ph)(x_k)\to(f\oplus\ph)(\ox)= f(\ow)+\ph(\ow-\ox)$$ and $\ph$ is continuous at $\ow-\ox$, $f(w_k)\to f(\ow)$. Thus
\begin{equation*}
x^*\in \partial f(\ow)\cap [-\partial\ph(\ow-\ox)].
\end{equation*}
The proof is now complete. $\h$

\section{Subdifferential Characterizations for Differentiability}

Let $f: X\to \oR$ be an extended-real-valued function with $\ox\in \mbox{\rm int dom}\,f$. We say that $f$ is \emph{Hadamard strictly differentiable at} $\ox$  if there exists $v\in X^*$ such that
\begin{equation*}
\lim_{x\to \ox, t\to 0^+}\dfrac{f(x+td)-f(x)-t\la v, d\ra}{t}=0,
\end{equation*}
where the convergence is uniform for $d$ in every compact subsets of $X$. The element $v$ is called the \emph{strict Hadamard derivative} of $f$ at $\ox$ denoted by $\nabla_Hf(\ox)$.

We can show that the Fr\'echet strict differentiability and the Hadamard strict differentiability are equivalent in finite dimensions.

We say that $\Hat\partial f(\cdot)$ is strongly continuous at $\ox$ if there exists an element $x^*\in X^*$ such that whenever $x_k\to \ox$ and $x^*_k\in \Hat\partial f(x_k)$, one has that $\|x^*_k-x^*\|\to 0$. It can be equivalently written as: there exists $x^*\in X^*$ such that for any $\varepsilon>0$, there exists $\delta>0$ such that whenever $\|x-\ox\|<\delta$ and $u^*\in \Hat\partial f(x)$, one has $\|u^*-x^*\|<\varepsilon$.

\begin{Theorem}\label{scon} Let $X$ be an Asplund space (see \cite{Mordukhovich_2006} for the definition) and let $f: X\to \oR$ be an extended-real-valued function with $\ox\in \mbox{\rm int dom}\,f$. Then the following are equivalent:\\[1ex]
{\rm\bf (i)} $f$ is locally Lipschitz continuous around $\ox$ and $\Hat\partial f(\cdot)$ is strongly continuous at $\ox$.\\
{\rm\bf (ii)} $f$ is Fr\'echet strictly differentiable at $\ox$.
\end{Theorem}
{\bf Proof.} Suppose that $f$ is locally Lipschitz continuous around $\ox$, $\Hat\partial f(\cdot)$ is continuous at $\ox$, and $f$ is not Fr\'echet strictly differentiable at $\ox$. Let $x^*$ be an element of the definition of strongly continuous of $\Hat\partial f(\cdot).$ Then, without loss of generality, we can assume that there exist $\gamma>0$ and sequences $x_k, y_k\to \ox$, $x_k\neq y_k$, such that
\begin{equation*}
\gamma\leq \lim_{k\to\infty}\dfrac{f(x_k)-f(y_k)-\la x^*, x_k-y_k\ra}{\|x_k-y_k\|}.
\end{equation*}
 By the mean value theorem \cite[Corollary 3.2]{Loewen} (with also holds in Asplund spaces; see \cite{Mordukhovich_2006}), there exist $c_k\to \ox$, $x^*_k\in \Hat\partial f(c_k)$ with
\begin{equation*}
f(x_k)-f(y_k)\leq \la x^*_k, x_k-y_k\ra +\|x_k-y_k\|^2.
\end{equation*}
Then
\begin{align*}
\gamma\leq\lim_{k\to \infty}\dfrac{f(x_k)-f(y_k)-\la x^*, x_k-y_k\ra}{\|x_k-y_k\|}&\leq \lim_{k\to \infty}\dfrac{\la x^*_k, x_k-y_k\ra +\|x_k-y_k\|^2-\la x^*, x_k-y_k\ra}{\|x_k-y_k\|}\\
&\leq \lim_{k\to \infty} (\|x_k-y_k\|+\|x^*_k-x^*\|)=0,
\end{align*}
which is a contradiction.

Now, we suppose that $f$ is Fr\'echet strictly differentiable at $\ox$ with $\nabla f(\ox)=x^*$. It is not hard to see that $f$ is locally Lipschitz continuous around $\ox$. Moreover, for any $\varepsilon>0$, there exists $\delta>0$ such that
$$\frac{f(x)-f(y)-\langle x^*,x-y\rangle }{\|x-y\|}\leq \left|\frac{f(x)-f(y)-\langle x^*,x-y\rangle }{\|x-y\|}\right|<\frac{\varepsilon}{2} $$
whenever $ \|x-\ox\|<\delta, \|y-\ox\|<\delta, x\ne y.$
So,
\begin{equation}\label{eq1thr6.1}
-\langle x^*,x-y\rangle\leq -f(x)+f(y)+\varepsilon\|x-y\|\ \text{whenever}\ \|x-\ox\|<\delta, \|y-\ox\|<\delta, x\ne y.
\end{equation}

Let $\delta'=\delta/2>0$ and let $y\in X$ such that $\|y-\ox\|<\delta'$ and  $u^*\in\Hat\partial f(y)$. We will show that $\|u^*-x^*\|\leq\varepsilon.$

It follows from $u^*\in\Hat\partial f(y)$ that there exists $\delta''<\delta'$ such that
\begin{equation}\label{eq2thr6.1}
\langle u^*,x-y\rangle\leq f(x)-f(y)+\frac{\varepsilon}{2} \|x-y\|\  \text{whenever}\ \|x-y\|<\delta''.
\end{equation}

If $x\in X$ such that $\|x-y\|<\delta'',$ then $\|x-\ox\|<\delta.$ It follows from (\ref{eq1thr6.1}) and (\ref{eq2thr6.1}) that
$$\langle u^*-x^*,x-y\rangle\leq\varepsilon\|x-y\|.$$

Therefore, $\|u^*-x^*\|\leq\varepsilon.$  $\h$

\begin{Corollary} Let $X$ be finite dimensional and let $f: X\to \oR$ be an extended-real-valued function with $\ox\in \mbox{\rm int dom}\,f$.
Then the following are equivalent:\\[1ex]
{\rm\bf(i)} $f$ is Hadamard strictly differentiable at $\ox$.\\
{\rm\bf (ii)} $f$ is Fr\'echet strictly differentiable $\ox$.\\
{\rm\bf (iii)} $f$ is locally Lipschitz continuous at $\ox$ and $\partial f(\ox)$ is a singleton.

Moreover, if $f$ is strictly differentiable on an open set $D$, then it is  continuously differentiable on this set.
\end{Corollary}
{\bf Proof.} The equivalence  {\bf(i)} $\Longleftrightarrow$ {\bf (ii)}  is well known and will be proved in Proposition \ref{eqHF} for the convenience of the reader. The implication {\bf (ii)} $\Longrightarrow$ {\bf (iii)} is trivial. In order to prove the implication {\bf(iii)} $\Longrightarrow$ {\bf(ii)}, by Theorem \ref{scon}, it suffices to show that
$\Hat\partial f(\cdot)$ is strongly continuous at $\ox$ under the assumption that $f$ is locally Lipschitz continuous at $\ox$ and $\partial f(\ox)$ is a singleton. Let $x^*$ be the only element of $\partial f(\ox).$ By contradiction, suppose that $\Hat\partial f(\cdot)$ is not strongly continuous at $\ox.$ Then there exists $\varepsilon_0>0$ and a sequence $\{x_k\}$ that converges to $\ox$ with $x_k^*\in\Hat\partial f(x_k)$ satisfying $\|x_k^*-x^*\|>\varepsilon_0$ for every $k$. Since $f$ is locally Lipschitz continuous at $\ox$, the sequence $\{x_k^*\}$ is bounded. Then there exists a subsequence $\{x_{k_l}^*\}$ of $\{x_k^*\}$ that converges to $y^*\in X^*.$   So $y^*\in\partial f(\ox)=\{x^*\}$ which yields a contradiction.  The last conclusion is trivial because the strong continuity of the Fr\'echet subdifferential mapping coincides with the continuity in this case. $\h$

For simplicity, we assume in what follows that $X$ is finite dimensional.

\begin{Proposition} In the setting of Theorem \ref{thr8} suppose that $X$ is finite dimensional. If $f$ is Fr\'echet strictly differentiable at $\ox$ or $\ph$ is Fr\'echet strictly differentiable at $0$, then $f\oplus\ph$ is Fr\'echet strictly differentiable at $\ox$.
\end{Proposition}
{\bf Proof.} Note that $\ph$ is convex and finite around $0$, so it is locally Lipschitz around $0$. Thus $f\oplus\ph$ is locally Lipschitz around $\ox$ and $\partial (f\oplus\ph)(\ox)$ is a singleton under the assumptions made, so it is Fr\'echet strictly differentiable at this point. $\h$

\begin{Proposition} Suppose that $X$ is finite dimensional, $\ph$ is Fr\'echet strictly differentiable, and $\mathcal{P}$ is inner semicompact at $\ox$. If $\mathcal{P}(\ox)$ is a singleton, then $f\oplus\ph$ is Fr\'echet strictly differentiable at $\ox$ and
\begin{equation*}
\nabla(f\oplus \ph)(\ox)=-\nabla \ph(\ow-\ox),
\end{equation*}
where $\ow\in \mathcal{P}(\ox)$.
\end{Proposition}
{\bf Proof.} Since $\ph$ is Fr\'echet strictly differentiable at $\ox$, it is locally Lipschitz continuous at this point, and so is $f\oplus\ph$. This implies that $\partial (f\oplus\ph)(\ox)$ is nonempty; see \cite[Corollary 2.25]{Mordukhovich_2006}. Then $\partial (f\oplus\ph)(\ox)$ is a singleton by Proposition \ref{lmt3}. $\h$

\begin{Example} {\rm Let $X=\R^n$ with the Euclidean norm and let $\ox\in X.$  Consider quadratic infimal convolution defined in  \eqref{quincon}.
Suppose that $f$ lower semicontinuous and is bounded below. Then we can show that $\mathcal{P}$ is inner semicompact at $\ox$. Suppose that $\mathcal{P}(\ox)=\{\ow\}$ (which holds if $f$ is convex). Then
$$\partial f_\alpha(\ox)\subset \partial f(\ow)\cap [-\nabla \ph(\ow-\ox)].$$
By Proposition \ref{LIP}, the function $f_\alpha(\ox)$ is locally Lipschitz continuous, so $\partial f_\alpha(\ox)$ is nonempty. It follows that
$$\partial f_\alpha(\ox)= [-\nabla \ph(\ow-\ox)]=2\alpha(\ox-\ow).$$
In fact, $f_\alpha$ is a $C^1$ function.
}

\end{Example}

\section*{Appendix: More on Subdifferential Characterizations for Differentiability}

In what follows we present some known results on subdifferential characterizations for differentiability; see, e.g., \cite{CL,bwang10}. Detailed proofs are given for the convenience of the reader.

\begin{Proposition}\label{pro6.7} Let $f:X\to \oR$ be an extended-real-valued convex function and let $\ox\in \mbox{\rm int dom}\,f.$ Then the following are equivalent:\\[1ex]
{\rm\bf(i)} $f$ is Hadamard strictly differentiable at $\ox$.\\
{\rm\bf (ii)} $f$ is locally Lipschitz continuous at $\ox$ and  G\^{a}taeux differentiable at $\ox$.\\
{\rm\bf(iii)} $f$ is locally Lipschitz continuous at $\ox$ $\partial f(\ox)$ is a singleton.
\end{Proposition}
{\bf Proof.} The implication {\bf (i)}$\Longrightarrow$ {\bf(ii)} follows from the definition and \cite[Proposition 2.2.1]{CL}. The implication {\bf (ii)}$\Longrightarrow$ {\bf(iii)} is obvious because if $f$  is convex and G\^{a}taeux differentiable at $\ox,$ then its subdifferential in the sense of convex analysis $\partial f(\ox)$ reduces to the  G\^{a}taeux derivative of $f$ at $\ox$ . The proof of the implication {\bf (iii)}$\Longrightarrow$ {\bf(i)}  can be found in \cite[Proposition 2.2.4]{CL} with the observation that   if $f$ is convex, then $\partial_Cf(\ox)=\partial f(\ox)$, where $\partial_Cf(\ox)$ denotes the Clarke subdifferential; see the definition in \cite{CL}. $\h$

\begin{Proposition} Let $f: X\to \oR$ be an extended-real-valued convex function  and let $\ox\in \mbox{\rm int dom}\,f.$ Then the following are equivalent:\\[1ex]
{\rm\bf(i)} $f$ is Fr\'echet strictly differentiable at $\ox$.\\
{\rm\bf (ii)} $f$ is locally Lipschitz continuous at $\ox$ and Fr\'echet differentiable at $\ox$.\\
{\rm\bf (iii)}  $f$ is locally Lipschitz continuous at $\ox,$  $\partial f(\ox)$ is a singleton, and $\partial f(\cdot)$ is strongly continuous at $\ox$.
\end{Proposition}
{\bf Proof.} The implication {\bf (i)}$\Longrightarrow$ {\bf(ii)} is obvious. If $f$ is locally Lipschitz continuous and Fr\'echet differentiable at $\ox,$  it is well-known that $\partial f(\ox)=\{\nabla f(\ox)\}$. Moreover, the subdifferential mapping is strongly continuous at $\ox$. Thus, the implication {\bf (ii)}$\Longrightarrow$ {\bf(iii)} holds. We now prove {\bf (iii)}$\Longrightarrow$ {\bf(i)}. Let $\partial f(\ox)=v^*.$ Since $\partial f(\cdot)$  is strongly continuous at $\ox,$ for any $\varepsilon>0$ there exists $\delta>0$ such that
$$\partial f(u)\subset \B(v^*,\varepsilon)\ \text{whenever}\ u\in\B(\ox,\delta).$$

We can choose $\delta>0$ sufficiently small such that   $f$ is Lipschitz continuous on $\B(\ox,\delta).$ Fix any $x,y\in \B(\ox; \delta)$ with $x\neq y$. By the subdifferential mean value theorem, there exist $u\in (x,y)$ and $w^*\in \partial f(u)$ such that
$$f(x)-f(y)=\la w^*, x-y\ra.$$
Then $\|w^*-v^*\|<\varepsilon$, and hence
\begin{equation*}
\big |\dfrac{f(x)-f(y)-\la v^*, x-y\ra}{\|x-y\|}\big|=\big |\dfrac{\la w^*-v^*, x-y\ra}{\|x-y\|}\big |\leq \|w^*-v^*\|<\varepsilon.
\end{equation*}
Thus, $f$ is Fr\'echet strictly differentiable at $\ox$. $\h$

The following known result follows from the fact that every bounded set in a
finite-dimensional space is contained in a compact set. Thus, the uniformity of
convergence with respect to bounded sets is implied by the uniformity
of convergence with respect to compact sets. We present here a direct proof.

\begin{Proposition} \label{eqHF} Suppose that $X$ is  finite dimensional.  Let $f: X\to \oR$ be an extended-real-valued  function  and let $\ox\in \mbox{\rm int dom}\,f.$ Then $f$ is Hadamard strictly differentiable at $\ox$ if and only if it is Fr\'echet strictly differentiable at $\ox$.
\end{Proposition}
{\bf Proof.} It is easy to see that the Fr\'echet strict differentiability implies the Hadamard strict differentiability. Let us prove the converse. By contradiction, suppose that $f$ is not Fr\'echet strictly differentiable at $\ox$. Then there exist $\varepsilon_0>0$ and sequences $x_k, y_k\to \ox$ with $x_k\neq y_k$ and
\begin{equation*}
\big |\dfrac{f(x_k)-f(y_k)-\la v, x_k-y_k\ra}{\|x_k-y_k\|}\big |\geq \varepsilon_0.
\end{equation*}
Let $d_k:=\dfrac{x_k-y_k}{\|x_k-y_k\|}$ and $t_k:=\|x_k-y_k\|$. Without loss of generality, suppose that $d_k\to d$ with $\|d\|=1$ as $k\to \infty$. Then
\begin{equation*}
\big |\dfrac{f(y_k+t_kd_k)-f(y_k)-\la v, t_kd_k\ra}{t_k}\big |\geq \varepsilon_0.
\end{equation*}
By \cite[Propsition 2.2.1]{CL}, $f$ is locally Lipschitz continuous around $\ox$ with Lipschitz constant $\ell$. Thus,
\begin{align*}
&\big |\dfrac{f(y_k+t_kd)-f(y_k)-\la v, t_kd\ra}{t_k}\big |\\
&=\big |\dfrac{f(y_k+t_kd)-f(y_k+t_kd_k)+f(y_k+t_kd_k)-\la v, t_kd_k\ra +\la v, t_kd_k\ra-f(y_k)-\la v, t_kd\ra}{t_k}\big |\\
&=\big |\dfrac{f(y_k+t_kd)-f(y_k+t_kd_k)+f(y_k+t_kd_k)-f(y_k)-\la v, t_kd_k\ra +\la v, t_kd_k\ra-\la v, t_kd\ra}{t_k}\big |\\
&\geq \big |\dfrac{f(y_k+t_kd_k)-f(y_k)-\la v, t_kd_k\ra}{t_k}\big |-\big |\dfrac{f(y_k+t_kd)-f(y_k+t_kd_k)}{t_k}\ra\big |-\big |\dfrac{\la v, t_kd_k\ra-\la v, t_kd\ra}{t_k}\big |\\
&\geq \varepsilon_0-\ell \|d_k-d\|-\|v\|\|d_k-d\|.
\end{align*}
It follows that
\begin{equation*}
\liminf_{k\to\infty}\big |\dfrac{f(y_k+t_kd)-f(y_k)-\la v, t_kd\ra}{t_k}\big |\geq \varepsilon_0,
\end{equation*}
which is a contradiction by \cite[Propsition 2.2.1]{CL}. $\h$

\begin{Example}{\rm Consider the function $f(x)=\|x\|$, $x\in \ell^1$. It is not hard to verify that $f$ is Hadamard strictly differentiable at every $x=(x_1, x_2, \ldots)\in \ell^1$, where $x_i\neq 0$ for every $i$, but it is not Fr\'echet strictly differentiable at that point.
}\end{Example}

The following corollaries can be derived easily.

\begin{Corollary} Let $X$ be finite dimensional and let $f: X\to \oR$ be an extended-real-valued convex function with $\ox\in \mbox{\rm int dom}\,f$. Then the following are equivalent:\\[1ex]
{\rm\bf(i)} $f$ is Hadamard strictly differentiable at $\ox$.\\
{\rm\bf (ii)} $f$ is  G\^{a}taeux differentiable at $\ox$.\\
{\rm\bf(iii)} $f$ is Fr\'echet strictly differentiable at $\ox$.\\
{\rm\bf(iv)} $f$ is Fr\'echet differentiable at $\ox$.\\
{\rm\bf(v)} $\partial f(\ox)$ is a singleton.
\end{Corollary}

\begin{Corollary} Let $X$ be finite dimensional and let $f: X\to \oR$ be an extended-real-valued convex function with $D:=\mbox{\rm int dom}\,f\neq \emptyset$. Then the following are equivalent:\\[1ex]
{\rm\bf(i)} $f$ is Hadamard strictly differentiable on $D$.\\
{\rm\bf (ii)} $f$ is Fr\'echet strictly differentiable on $D$.\\
{\rm\bf(iii)} $f$ is continuously differentiable on $D$.\\
{\rm\bf(iv)} $\partial f(x)$ is a singleton for all $x\in D$.
\end{Corollary}

\end{document}